\documentclass[prd,eqsecnum,showpacs,nofootinbib]{revtex4}
\usepackage{bm}
\usepackage{amsmath}
\usepackage{amssymb}
\usepackage{graphicx}
\usepackage{epstopdf}

\begin{document}
\title{Simplicial Ricci Flow: An Example of a Neck Pinch Singularity in 3D}\author{Paul M. Alsing$^{1}$, Warner A. Miller$^{2,3}$\footnote{Corresponding author: wam@fau.edu}, Matthew Corne$^1$, David Gu$^{4}$,\\  Seth Lloyd$^5$, Shannon Ray$^3$  \& Shing-Tung Yau$^2$}
\affiliation{$^1$ Air Force Research Laboratory, Information Directorate, Rome, NY 13441\\
$^2$ Department of Mathematics, Harvard University, Cambridge MA 02138\\
$^3$ Department of Physics, Florida Atlantic University, Boca Raton, FL 33431\\
$^4$ Department of Computer Science, Stony Brook University, Stony Brook, NY 11794\\
$^5$ MIT Department of Mechanical Engineering, MIT 3-160,  Cambridge, MA 02139}
\date{Submitted: June 25, 2013}
\begin{abstract}
We examine a Type-1 neck pinch singularity in simplicial Ricci flow (SRF) for an axisymmetric piecewise flat 3-dimensional geometry with topology $S^3$.  SRF was recently introduced as an unstructured mesh formulation of Hamilton's Ricci flow (RF).  It describes the RF of a piecewise-flat simplicial geometry.  In this paper, we apply the SRF equations to a representative double-lobed axisymmetric piecewise flat geometry with mirror symmetry at the neck similar to the geometry studied by Angenent and Knopf (A-K).   We choose a specific radial profile and compare the SRF equations with the corresponding finite-difference solution of the continuum A-K RF equations. The piecewise-flat 3-geometries considered here are built of isosceles-triangle-based frustum blocks. The axial symmetry of this model allows us to use frustum blocks instead of tetrahedra.  The $S^2$ cross-sectional geometries in our model are regular icosahedra. We demonstrate that, under a suitably-pinched initial geometry, the SRF equations for this relatively low-resolution discrete geometry yield the canonical Type-1 neck pinch singularity found in  the corresponding continuum solution.   We adaptively remesh during the evolution to keep the circumcentric dual lattice well-centered. Without such remeshing, we cannot evolve the discrete geometry to neck pinch. We conclude with a discussion of future generalizations and tests of this SRF model.
\end{abstract}
\pacs{04.60.Nc,02.40.Hw, 02.40.Ma,02.40.Ky}
\maketitle

\section{Exploring Simplicial Ricci Flow in 3D}

Hamilton's Ricci flow (RF) has yielded new insights into pure and applied mathematics as well as engineering fields \cite{Hamilton:1982,Cao:2003,Chow:2004,Chow:2006,Chow:2007}.  Here the time evolution of the metric is proportional to the Ricci tensor,
\begin{equation}\label{eq:RF}
\dot g_{ab} = -2\, R_{ab},
\end{equation}
and yields a forced diffusion equation for the curvature; i.e., the scalar curvature evolves as
\begin{equation}\label{eq:fdeqn}
\dot R = \triangle R + 2 R^2.
\end{equation}
The bulk of the applications of this curvature flow technique have been limited to the numerical evolution of piecewise-flat simplicial 2--surfaces \cite{Gu:2012}.  It is well established that,  ordinarily, a geometry with complex topology is most naturally represented in a coordinate-free way by unstructured meshes.  e.g. finite volume \cite{Peiro:2005}, finite element \cite{Humphries:1997}, in general relativity by  Regge calculus \cite{Regge:1961,Gentle:1998} and for electrodynamics  by discrete exterior calculus  \cite{Hirani:2005}.  While the utility of piecewise--flat simplicial geometries in analyzing the RF of 2--dimensional geometries is well established and proven to be effective \cite{Gu:2012,Chow:2003}, one expects a wealth of exciting new applications in 3 and higher dimensions.  Here we explore the utility of RF in three and higher dimensions.  However, we first note that the applications of discrete RF in two dimensions arise from its diffusive curvature properties and from the uniformization theorem, that every simply connected Riemann surface evolves under RF to one of three constant curvature surfaces --- a sphere, a Euclidean plane or a hyperbolic plane.  RF on surfaces is perhaps the only general method to engineer a metric for a surface given only its curvature \cite{Gu:2012}.  In three dimensions the uniformization theorem yields to the geometrization conjecture of Thurston suggesting that all Riemannian 3-manifolds have a similar,  but richer,  classification into a connected sum of one or more of eight canonical geometries \cite{Besson:2007}.  The diffusive curvature flow in 3 and higher dimensions together with this classification can provide a richer taxonomy than its 2-dimensional counterpart.  We believe this more refined taxonomy will prove useful in network classification.  Diffusive curvature flow can provide noise reduction in higher dimensional manifolds, and in this direction we are currently exploring a coupling RF with persistent homology \cite{Corne:2014,Mischaikow:2013}.  

A discrete RF approach for three and higher dimensions, referred to as Simplicial Ricci Flow (SRF), has been introduced recently and is founded on Regge calculus \cite{Miller:2013,AMM:2011,McDonald:2012}, as well as complementary work in this direction by \cite{Glickenstein:2011a,Glickenstein:2011,G:2005,Ge:2013,Forman:2003,Glickenstein:2005,GuSaucan:2013}.  The Regge-Ricci flow (RRF) equations of SRF are similar to their continuum counterpart.  They are naturally defined on a $d$-dimensional simplicial geometry as a proportionality between the time rate of change of the circumcentric dual edges, $\lambda_i$, and the  simplicial Ricci tensor associated to these dual edges,
\begin{equation}\label{eq:ref}
\dot \lambda_i = -2\ Ric_\lambda.
\end{equation}
It is the aim of this paper to explore the behavior of these RRF equations in 3-dimensions for a geometry with axial symmetry, and to examine the development of a Type-1 neck pinch singularity.  We use as a foundation of this work the analysis of Angenent and Knopf on the Type-1 singularity analysis of the continuum RF equations \cite{Knopf:2004}.   They carefully analyzed a class of axisymmetric dumbbell-shaped geometries with mirror symmetry about the plane of the neck as illustrated in the top of Fig.~ \ref{fig:dumbbell}.  The metrics these researchers evolve under RF are commonly referred to as  warped product metrics on $I \times S^2$,
\begin{eqnarray}
\label{eq:wpm}
g & = &\underbrace{\varphi(z)^2 dz^2}_{da^2} + \rho(z)^2 g_{can}\\
 \label{eq:asmetric}
   & = & da^2 + \rho(a)^2 g_{can}.
\end{eqnarray}
Here, $I \in \mathbb{R}$ is an open interval,
\begin{equation}\label{eq:gcan}
g_{can} = d\theta^2 + \sin^2\theta d\phi^2,
\end{equation}
 is the metric of the unit 2-sphere,
\begin{equation}
a(z) = \int_{0}^{z} \, \varphi(z) dz,
\end{equation}
is the geodesic axial distance away from the waist, and $\rho(a)$ is the radial profile of the mirror-symmetric geometry, i.e. $s=\rho(a)$ is the radius of the cross-sectional 2-sphere at an axial distance $a$ away from the waist.  Angenent and Knopf proved  that the RF evolution for such a geometry has the following properties:
\begin{enumerate}
\item If the scalar curvature is everywhere positive, $R\ge0$,  then the radius of the waist ($s_{min}=\rho(0)$)  is bounded,
\begin{equation}
(T-t) \le s_{min}^2 \le 2(T-t),
\end{equation}
where $T$ is the finite time that a neck pinch occurs.
\item As a consequence, the neck pinch singularity occurs at or before $T = s_{min}^2$.
\item The height of the two lobes are bounded from below and under suitable conditions, the neck will pinch off before the lobes will collapse.
\item The neck approaches a cylindrical-type singularity.
\end{enumerate}
We demonstrate in this paper that the SRF equations, for a sufficiently pinched radial profile, will produce a neck pinch singularity in finite time.  Furthermore, we show that the results agree with a finite-difference solution of the continuum RF equations for the same profile. The discrete model here is a very coarse approximation to the dumbbell geometry (e.g. the $S^2$ cross sections are modeled by icosahedra, and adjacent faces of the icosahedra are connected to each other via frustum blocks).  However, this work represents the first non-trivial numerical solution of the SRF equations, and it is the goal of this paper to demonstrate for the first time neck pinch behavior in SRF.

We can very well believe that the set of RRF equations will have an equally rich spectrum of application as does its 2-dimensional counterpart known as combinatorial RF \cite{Chow:2003}.  We therefore are motivated to explore the discrete RF in higher dimensions so that it can be used in the analysis of topology and geometry, both numerically and analytically to bound Ricci curvature  in discrete geometries and to analyze and handle higher--dimensional RF singularities \cite{LinYau:2010,Knopf:2009}.  The  topological taxonomy afforded by RF is richer when transitioning from 2 to 3--dimensions.  In particular,  the uniformization theorem says that any  2--geometry will evolve under RF to a constant curvature sphere, plane or hyperboloid, while in 3--dimensions the curvature and surface will diffuse into a connected sum of prime manifolds \cite{Thurston:1997}.

\section{A Simplicial Approximation of a Angenent-Knopf Neck pinch-Type  Model: Initial Value Data and Lattice Structure}\label{sec:2}

\begin{figure}
\includegraphics[width=7.5cm]{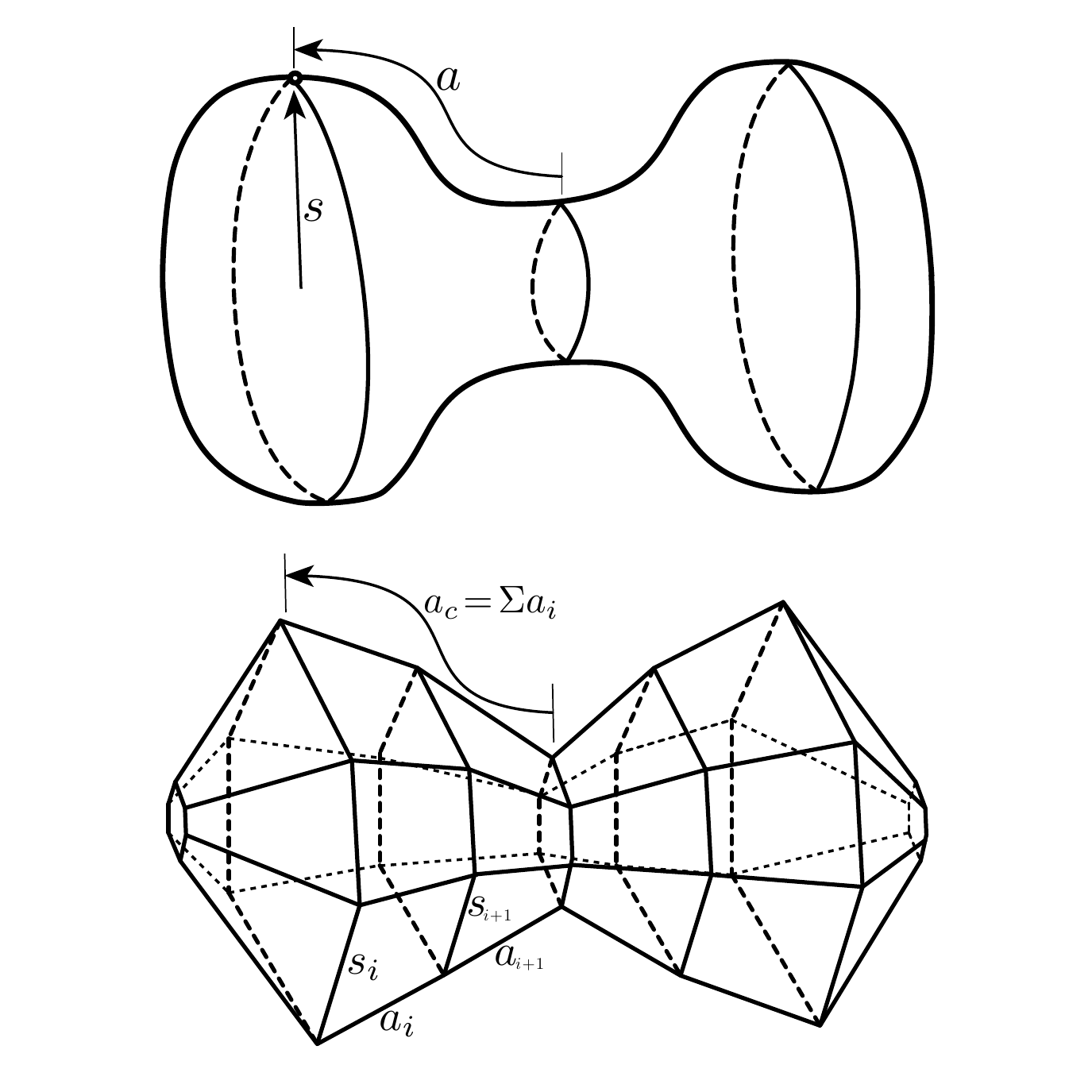}
\caption{An illustration of a 2-D dumbbell geometry analogous to the 3-D dumbbell we will analyze in this manuscript.  The top illustration shows an embedding surface of a continuum dumbbell, while the bottom figure shows a simplicial representation of this geometry using trapezoids of three shapes and two small hexagonal end-caps at the poles.  The topology of these models is $S^2$.  Just as this 2-D piecewise-flat dumbbell geometry is characterized by the sum of axial edges, $a_c = a^{(coord)}_i$, measuring the proper distance a hexagonal cross section (of edge length $s_i$)  is away from the waist;  the 3-D dumbbell is parameterized with the same set of edge lengths.  However, in the case of the 3-D dumbbell the surface is tiled with regular triangle frustum blocks and two icosahedral  end caps as shown Fig.~\ref{Fig:9icosa}.}
\label{fig:dumbbell}
\end{figure}

\begin{figure}
\includegraphics[width=10cm]{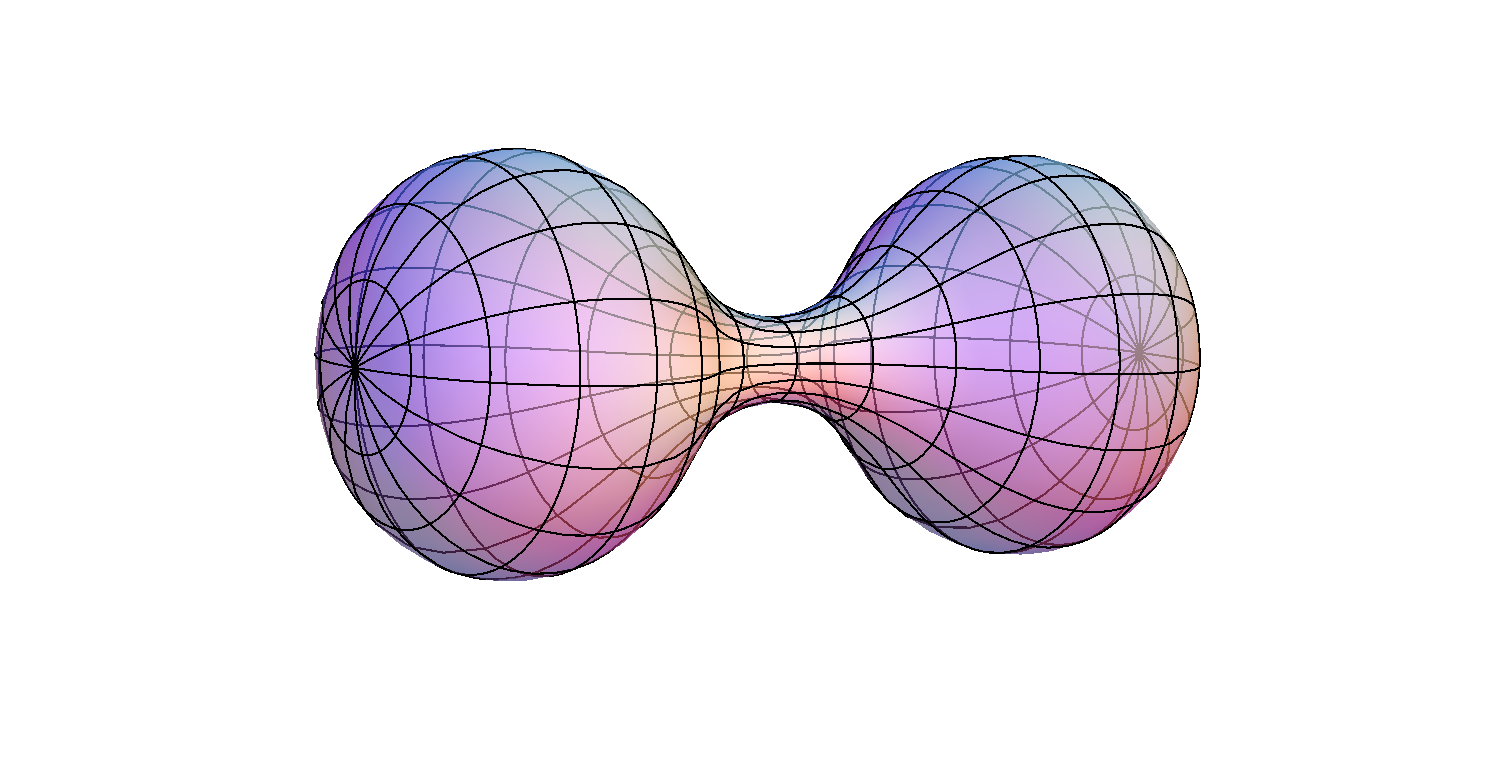}
\caption{We illustrate here an embedding of the initial 3-D hypersurface representing the dumbbell geometry at $t=0$.  This is the initial-value data  we examine in this manuscript.  It is determined by the radial profile in Eq.~\ref{Eq:ivd}. Here, one of the azimuthal angles have been suppressed so that it can be embedded in $R^3$.  The cross sections perpendicular to its symmetry axis are 2-spheres of various radii, $\rho$, and are given by the radial profile equation. In the case considered in this paper, the waist has the smallest radius, $\rho_w(t=0):=\rho_0=10$ at $a_c=a_w=0$, and the lobes have a maximal radius, $\rho_m(t=0) \approx 48.9$ at $a\approx \pm 86.0$.}
\label{Fig:cos4}
\end{figure}

For the purpose of examining Type-1 neck pinch behavior of the SRF equations introduced in \cite{Miller:2013}, we have adopted for our simulations the qualitative features of the Angenent and Knopf initial data \cite{Knopf:2004}, as shown in Fig.~\ref{fig:dumbbell}.  The cross sections of this geometry in planes perpendicular to the symmetry axis are 2-sphere surfaces.  As discussed in the last section, the surface and metric can be parameterized by two coordinates, $a$ and $s$.  Here a given point on the surface is identified by its proper ``axial" distance of $a$ from the throat, and the radius, $s=\rho(a)$, of the cross-sectional sphere where the point lies, Eqs.~\ref{eq:asmetric} and \ref{eq:gcan}.
The initial data is determined by a radial profile function for the dumbbell geometry, and amounts to specifying a function relating the cylindrical radius of the dumbbell, $s$,  to a scaled proper axial distance along the dumbbell away from the neck, $a\in\{-\pi/2,\pi/2\}$,
\begin{equation}
s = \rho(a).
\end{equation}
If the dumbbell geometry has no neck, and is just a sphere of radius, $R_0$, this initial value data (ivd) radial profile function is simply the cylindrical coordinate radius,
\begin{equation}
s = \rho_{sphere}(a)=R_0\, cos(a).
\end{equation}
In order to cleanly prove their theorems on the evolution of such geometries,  Angenent and Knopf introduced a parabolic waist for the purpose to aid in the mathematical analysis of the neck singularity,
\begin{equation}
r=\rho_{_{AK}} = \left\{
\begin{array}{l l}
R_0 \cos{(a)}  & |a|\ge \frac{\pi}{4} \\
R_0 \sqrt{A + B a^2} & |a| < \frac{\pi}{4}
\end{array}
\right. .
\end{equation}
The constant $A$  controls the degree of neck pinching in the ivd, while the constant $B$ is chosen so as to ensure continuity in the radial profile function at $a=\pm \pi/4$.  To better serve our purposes in this paper of numerically comparing the continuum RF with the RRF equations in SRF, we analyze a geometry with a smoother radial profile,
\begin{equation}\label{Eq:ivd}
\boxed{\rho(a) = R_0 \left(\cos{\left(a/R_0\right)} + (1-\rho_0) \cos^4{\left( a/R_0 \right)}\right), \ \  \forall\, a\in R_0 \left\{- \frac{\pi}{2}, \frac{\pi}{2}\right\}. }
\end{equation}
We examine here the single case where $\rho_0=0.1$ and $R_0=100$. A 2-D embedding diagram for this 3-D geometry (suppressing one of the azimuthal angles) is shown in Fig.~\ref{Fig:cos4}.

\begin{figure}
 \includegraphics[width=14.5cm]{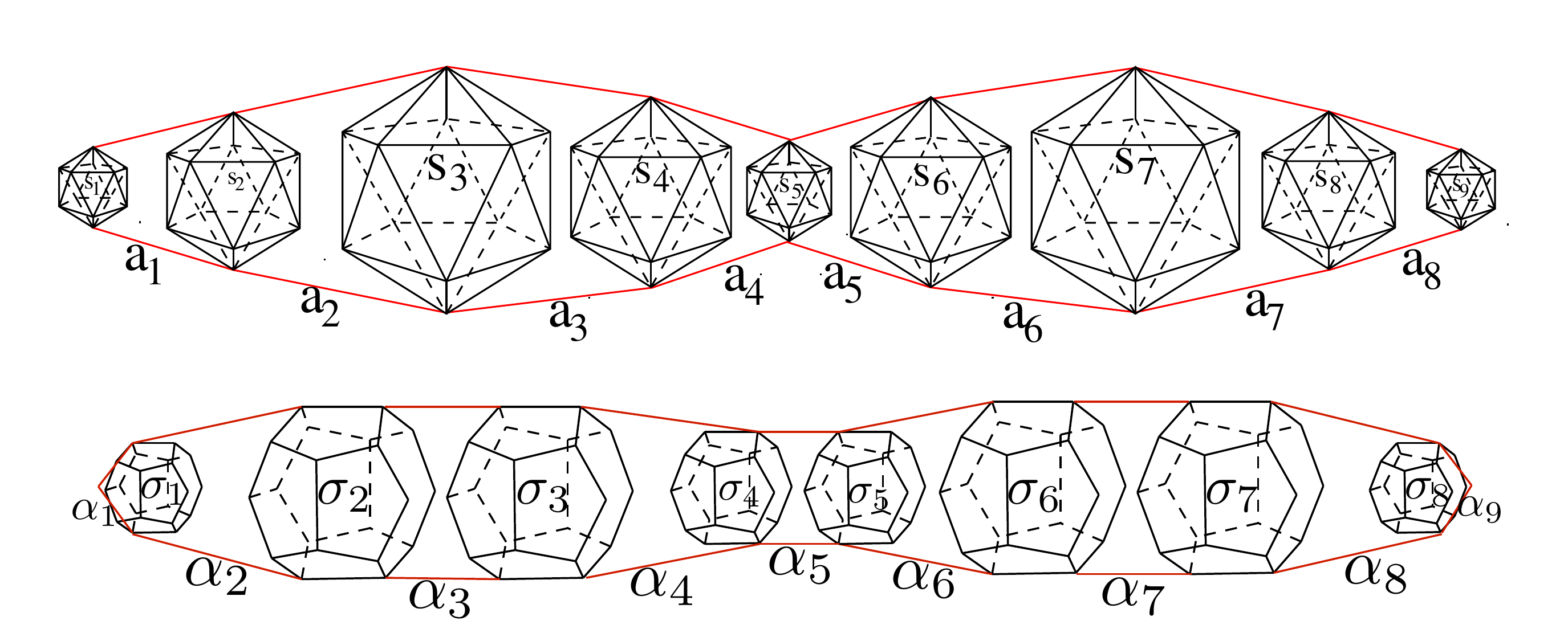}
\caption{{\em Top Figure:} The illustration at the top is our icosahedral dumbbell model with $n_s=9$ icosahedra.  The geometry of each of these icosahedra is determined by a single parameter, its edge length, $s_i$.  There are nine of these edge lengths, $\{s_1,s_2, \ldots, s_{n_s}\}$, in this model.  Adjacent icosahedra are connected with twelve equal axial edges of length $a_i$.  There are $n_s-1=8$ axial edges. For clarity we only show two of the twelve axial edges $a_i$ sandwiched between the $i$'th and $(i+1)$'st icosahedron. This yields a piecewise-flat geometry built of  $20 (n_s-1)$ frustum blocks (as described in Appendix~\ref{app:1}) and two icosahedral end caps.  {\em Bottom Figure:}  The bottom figure illustrates its circumcentric dual dodecahedral frustum lattice as described in the text.}
\label{Fig:9icosa}
\end{figure}

We approximate this double-lobed geometry (as shown in Fig.~\ref{Fig:cos4})  by first identifying a finite number, ($N_{s}$), of spherical cross sections along its symmetry axis.  To help visualize our lattice geometry, we illustrate this 3-D dumbbell model for the case of $n_s=9$ as shown in Fig.~\ref{Fig:9icosa}.  However,  in this paper we choose $n_s=45$ for our numerical simulation.  Each 2-sphere cross section is approximated by a regular icosahedron, of edge length $s_i$, where $i=1,2,\ldots, 45$.  The icosahedron with index $i=(n_s+1)/2=23$ is placed at the waist while the remaining 44 icosahedra are paired off so that each pair ($i=23\pm j$ for $j=1,2,\ldots, 22$) is placed at an ever increasing, but equal distance to the left and right of the waist.  We do not enforce mirror symmetry in our SRF equations; however, our initial value profile at $t=0$ is mirror symmetric.  The two end-cap icosahedra with index $i=1$ and $i=45$ are placed inside the two extreme ends of the dumbbell geometry.  In this way they initially have a nonzero radius, $s_{1}$ and $s_{45} > 0$.

We connect the vertices of each pair of adjacent icosahedra, i.e. icosahedra of index $i$ and $i+1$, using twelve equal axial edges of length, $a_i$.  We do not allow any twist of one icosahedron with respect to another; therefore, each regular triangle face of the $i$'th icosahedron, when connected by three axial edges, $a_i$, to its adjacent triangle on the $(i+1)$'st icosahedron will form a regular triangle frustum block as illustrated in Fig.~\ref{Fig:frustum}. There are twenty such frustum blocks sandwiched between any two adjacent icosahedra. The geometry of this piecewise-flat icosahedral frustum model  contains  $20(n_s-1)=880$ equilateral-triangle-based frustum blocks and two icosahedra. The lattice is rigid due to its axial symmetry and because we do not allow any twisting of the frustum blocks.  The geometry is completely determined by the $n_s=45$ icosahedral edge lengths, $s_i$ \, and the $(n_s-1)=44$ axial lengths, $a_i$.  Each evolution step requires the solution of $(2n_s-1)=89$ locally-coupled nonlinear algebraic first order SRF equations.  These equations are described in some detail in Sec.~\ref{sec:srfeqns}.

In order to compare our geometry with the continuum we set the circumradius of each icosahedron equal to the radius of sphere, $r_i$,
\begin{equation}\label{eq:s2r}
 r_i =\frac{\sqrt{10+2\sqrt{5}}}{4} \, s_i,
\end{equation}
and the proper distance from the waist $a$  in the continuum is related to the sum of the axial edges,
\begin{equation}\label{eq:a2ac}
a^{(coord)}_i =  \sum_{j = 1}^i a_i - \underbrace{\sum_{j=1}^{\frac{n_s-1}{2}} a_i}_{a_w}.
\end{equation}

The RRF equations will use extensively the elements of the circumcentric-dual lattice. The original icosahedron-based lattice described above is constructed by regular triangle-based frustum blocks and two end cap icosahedra.  Recall that in this icosahedral frustum lattice, there are $n_s$ icosahedral edge lengths, $s_i$,  and $(n_s-1)$ axial edges, $a_i$.  The dual lattice shown in Fig.~\ref{Fig:9icosa} is an axisymmetric geometry consisting of $(n_s-1)$ cross-sectional spheres approximated by regular dodecahedra with edges of length, $\sigma_i$.  There is one dodecahedron of the dual lattice sandwiched between every pair of adjacent icosahedra.  The 20 vertices of each dodecahedron are the 20 circumcenters of the 20 triangle frustum blocks sandwiched between two adjacent icosahedra.  Any two adjacent dodecahedra (e.g. the $i$'th and the $(i+1)$'st) in the dual lattice are connected to each other by the 20 dual axial edges, $\alpha_i$.  At each of the two end caps of the dumbbell geometry the 20 axial edges meet at a vertex. As mentioned, the vertices of the circumcentric dual lattice are located at the circumcenters of the triangular frustum blocks; however, at the two end caps the two dual lattice vertices are located at the circumcenter of the icosahedral end caps. The edges of the dual lattice connect adjacent circumcenters.  The geometry of the dual dodecahedral lattice is composed of $12(n_s-2)$ regular pentagonal-based frustum blocks and 12 tetrahedra at each end cap. Ordinarily, there are may more circumcentric dual edges than original edges; however, the discrete 3-dimensional axisymmetric geometry we are considering in this manuscript has an equal number.  We take advantage of this and show in Appendix~\ref{app:B} that the dual edge RRF equations are equivalent to the original RRF equations.

\begin{figure}
\includegraphics[width=7cm]{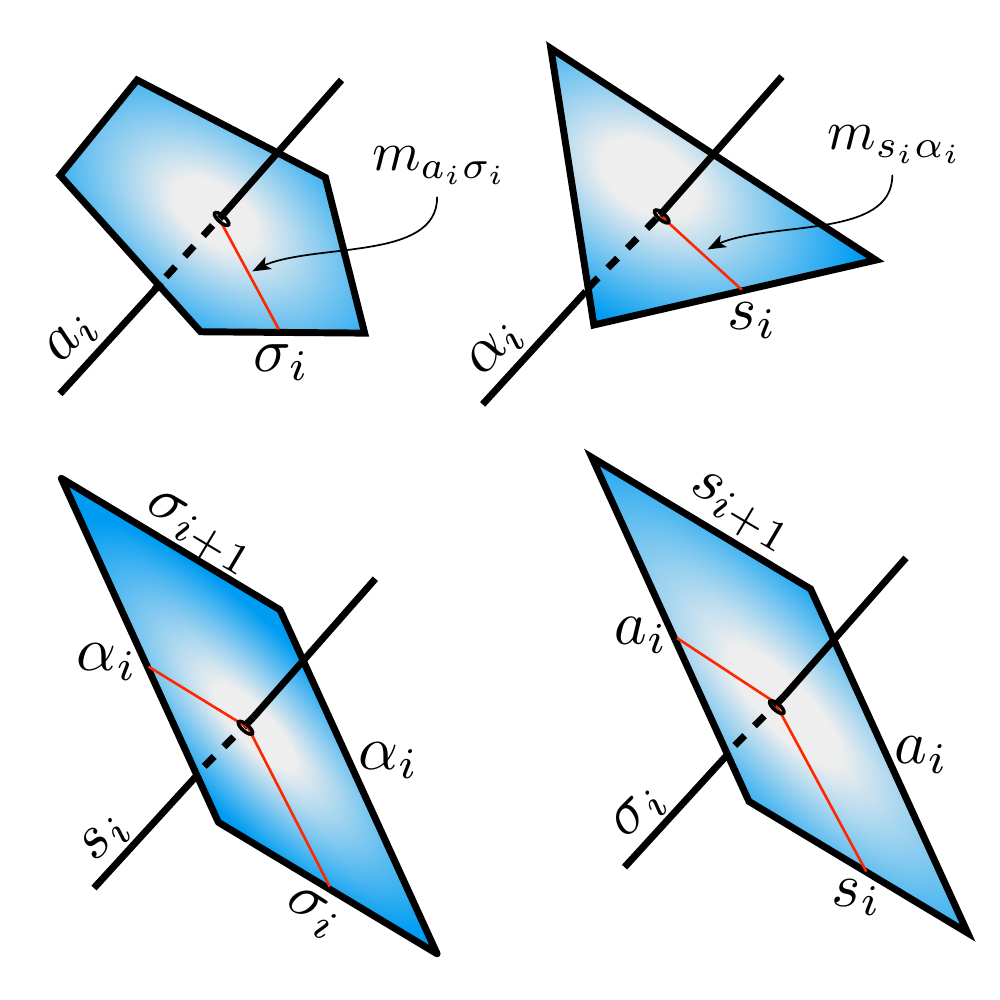}
\caption{We illustrate here the dual areas (shaded in blue) for each of the four types of edges in our lattice geometry.  We also identify a couple of representative moment arms (thin red lines) connecting the center of the edge, $a_i$ and $s_i$ to the perpendicular bisector of the dual edge, $\sigma_i$ and $\alpha_i$; respectively.  The pentagonal-shaped dual area $a^*_i$  to axial edge $a_i$ is used in part to define the sectional curvature associated of edge $a_i$ and is used in the construction of the Ricci tensor.  This is also true for the trapezoidal dual area $s^*_i$ to icosahedral edge $s_i$.}
\label{Fig:dualareas}
\end{figure}

Each edge of the icosahedral frustum lattice is dual to a polyhedron in the dodecahedral frustum lattice, and visa versa.  The dual areas are used to construct the RRF equations. In particular each  axial edge, $a_i$, is dual to a pentagonal face of the $i$'th dodecahedron. We refer to this dual pentagonal face by, $a^*_i$. Each icosahedral edge, $s_i$, is dual to a trapezoidal face $s^*_i$ of the $i$'th pentagonal-based frustum block.  The edges of $s^*_i$ are $\sigma_i$, $\sigma_{i+1}$ and $\alpha_i$.  Similarly, each dodecahedron edge, $\sigma_i$, is dual to a trapezoidal face $\sigma^*_i$  of the $i$'th regular triangle-based frustum block. The edges of $\sigma^*_i$ are edges $s_i$, $s_{i+1}$ and $a_i$.  Finally, each dual axial edge, $\alpha_i$, is dual to a regular triangle face of the $i$'th icosahedron, $\alpha^*_i$.  These dual areas and their associated edges are illustrated in Fig.~\ref{Fig:dualareas}

\section{The Continuum Ricci Flow Equations for a Warped Product Metric}\label{sec:continuum}

 In this section, we briefly highlight the continuum equations for the warped-product metric as introduced by Angenent and Knopf \cite{Knopf:2004} and Simon \cite{Simon:2000}. We do this so that we can more precisely compare the numerical solution of these continuum RF equations governing the dynamics of the specific radial profile (Eq.~\ref{Eq:ivd}) to that of the evolution of the corresponding SRF equations for our axisymmetric piecewise-flat lattice geometry. The numerical comparison will be discussed in Sec.~\ref{sec:numerics}, here we briefly outline the continuum equations.

There are only two distinct mixed Ricci tensors for the warped product metric given by Eq.~\ref{eq:asmetric},
\begin{eqnarray}
Rc^a{}_{a} & =&  -2 \left(\frac{\rho''}{\rho}\right),\\
Rc^\theta{}_\theta & = & Rc^\phi{}_\phi = \frac{1}{\rho^2} - \frac{\rho''}{\rho} -\left(\frac{\rho'}{\rho}\right)^2,
\end{eqnarray}
where  $\rho=\rho(a)$ and the primes are partial derivates with respect to the proper distance $a$. The RF equations for this warped-product metric in mixed form govern the dynamics of the radial profile by a single partial differential equation for $\rho(a)$,
\begin{equation}\label{eq:rhoeqn}
\frac{\dot \rho}{\rho} = - \frac{1}{\rho^2}+ \frac{\rho''}{\rho}  + \left(\frac{\rho'}{\rho}\right)^2.
\end{equation}
where the dots are time derivatives and the primes are partial derivatives with respect to $a$.  The time evolution of $\varphi$ can be recovered solely from the radial profile $\rho$ and its second derivatives, $\rho''$,
\begin{equation}
\frac{\dot \varphi}{\varphi} = - 2\,  \frac{\rho''}{\rho}.
\end{equation}
 We numerically solve Eq.~\ref{eq:rhoeqn} using the initial radial profile given by Eq.~\ref{Eq:ivd}.  This allows us in Sec.~\ref{sec:numerics}  to compare the solution of our SRF equations  with the continuum solution.

\section{The SRF equations for the Icosahedral Frustum Model}\label{sec:srfeqns}

For the axisymmetric model we are analyzing, we show in Appendix~\ref{app:B} that the RRF equations for the axial edges, $a_i$, and the icosahedral edges, $s_i$, are equivalent to the dual-edge RRF equations associated with dual edges $\sigma_i$ of the dodecahedra and dual axial edges $\alpha_i$ (Fig.~\ref{Fig:9icosa}).  In particular, we show that
\begin{equation}\label{eq:drrf-rrf}
\underbrace{
\left\{
\begin{array}{l}
\frac{\dot \sigma_i}{\sigma_i} = -Rc_{\sigma_i}\\
\frac{\dot \alpha_i}{\alpha_i} = -Rc_{\alpha_i}
\end{array}
\right\}
}_{dual-edge\ RRF\ equations}
\Longleftrightarrow \ \ \
\underbrace{
\left\{
\begin{array}{l}
\sum_{\lambda_j \in s^*_i}  \frac{\dot \lambda_j}{\lambda_j} \left(\frac{V_{s_i \lambda_j}}{V_{s_i}}\right) = -Rc_{s_i}\\
\sum_{\lambda_j \in a^*_i}  \frac{\dot \lambda_j}{\lambda_j} \left(\frac{V_{a_i \lambda_j}}{V_{a_i}}\right) = - Rc_{a_i}
\end{array}
\right\}
}_{simplicial-edge\ RRF\ equations},
\end{equation}
where the $\lambda_i$'s in the sum are the edges of the dual lattice in the boundary of the polygon, $a^*_i$ or $s^*_i$ dual to the edge $a_i$ or $s_i$, respectively (as shown in Fig.~\ref{Fig:dualareas}). The dual-edge RRF equations are significantly simpler.  For this reason, we have chosen to solve the dual-edge equations for this model.  Our analysis relies heavily on a recent paper by some of us, and in particular on Sec. 3 and Sec. 5 of \cite{Miller:2013}. The notation here corresponds to the notation used in the recent SRF manuscript \cite{Miller:2013}.

\subsection{The $\alpha_i$ dual axial equations}

There is one axial dual-edge equation associated each of the $n_s$ dual axial edges $\alpha_i$.  Each of these equations can be expressed locally in terms of the four edges, $a_i$, $a_{i+1}$, $s_i$ and $s_{i+1}$ of the icosahedral frustum model and their time derivatives,
\begin{equation}\label{eq:alpha}
\frac{\dot \alpha_i}{\alpha_i} = -Rc_{\alpha_i}.
\end{equation}
The axial edge, $\alpha_i$, reaches from the circumcenter of one triangle frustum to the adjacent triangle frustum sharing an equilateral triangle of the $i$'th icosahedron and is given by the sum of Eqs.~\ref{eq:theta_{s_i}} and \ref{eq:theta_{s_{i+1}}}, except for the two dual axial edges on the end caps where they terminate at the circumcenter of the end cap icosahedron. In particular, we find
\begin{equation}
\alpha_i = \left\{
\begin{array}{ll}
 \frac{\sqrt{3}}{12}\left(3+\sqrt{5}\right) s_1 + \frac{\sqrt{3}}{6} \left( \frac{3 a_1^2-2 s_1 (s_1-s_2)}{\sqrt{3 a_1^2-(s_1-s_2)^2}}\right), & i=1, \\
\frac{\sqrt{3}}{6} \left( \frac{3 a_i^2+2 s_{i+1} (s_i-s_{i+1})}{\sqrt{3 a_i^2-(s_i-s_{i+1})^2}}\right)+ \frac{\sqrt{3}}{6} \left( \frac{3 a_{i+1}^2-2 s_{i+1} (s_{i+1}-s_{i+2})}{\sqrt{3 a_{i+1}^2-(s_{i+1}-s_{i+2})^2}}\right), & i\in\{2,3,\ldots,N_{s-1}\},\\
\frac{\sqrt{3}}{12}\left(3+\sqrt{5}\right) s_{n_s} + \frac{\sqrt{3}}{6} \left( \frac{3 a_{n_s}^2-2 s_{n_s} (s_{n_s}-s_{n_s-1})}{\sqrt{3 a_{n_s}^2-(s_{n_s}-s_{n_s-1})^2}}\right), & i=n_s,
\end{array}
\right.
\end{equation}
together with the time derivatives,
\begin{equation}
\dot \alpha_i = \frac{\partial \alpha_i}{\partial a_i} \, \dot a_i + \frac{\partial \alpha_i}{\partial a_{i+1}} \, \dot a_{i+1} +  \frac{\partial \alpha_i}{\partial s_i} \, \dot s_i +   \frac{\partial \alpha_i}{\partial s_{i+1}} \, \dot s_{i+1}.
\end{equation}
These quantities are used to construct the left-hand side of Eq.~\ref{eq:alpha}.

The right-hand side of Eq.~\ref{eq:alpha} are the Ricci tensors associated with the dual axial edges,
\begin{equation}
Rc_{\alpha_i} = \sum_{\ell_j \in \alpha^*_i} Rc^{(hyb)}_{\ell_j} \left( \frac{V_{\alpha_i \ell_j}}{V_{\alpha_i}} \right) = 3  Rc^{(hyb)}_{s_i} \left( \frac{1}{3} \right) = 2 Rm_{s_i} = 2\, \frac{\epsilon_{s_i}}{s^*_i}.
\end{equation}
Here, the deficit angle $\epsilon_{s_i}$ is associated with the icosahedral edge, $s_i$.  Edge $s_i$  is the edge common to four triangle frustum blocks, except for the bounding icosahedral caps where it is the edge common to two frustum blocks and one end cap icosahedron; we find
\begin{equation}\label{eq:defs}
\epsilon_{s_i} = \left\{
\begin{array}{ll}
2\pi - \theta_{icosa} - 2\ \theta{s_1}, & i=1\\
2\pi - 4\,\theta{s_{i}}, & i\in\{2,3,\ldots,N_{s-1}\},\\
\theta_{icosa} + 2\,\theta{s_{n_s-1}}. & i=n_s,
\end{array}
\right.
\end{equation}
The dihedral angles ($\theta$)  are give in Appendix~\ref{app:1}. Care must be taken when calculating the dual areas, $s^*_i$, as shown in the lower-right of Fig.~\ref{Fig:dualareas}.  This dual area is given by moment arms and the three distinct dual edges in the boundary of $s^*_i$. It is simply the sum of four triangle areas, two of which are always equal,
\begin{equation}
s^*_i = \left\{
\begin{array}{ll}
m_{s_1\alpha_1}\, \alpha_1 + \frac{1}{2} m_{s_1 \sigma_1}, & i=1\\
m_{s_i\alpha_i}\, \alpha_i + \frac{1}{2} m_{s_i \sigma_i}\,  \sigma_i + \frac{1}{2} m_{s_i \sigma_{i+1}}\, \sigma_{i+1}, & i\in\{2,3,\ldots,N_{s-1}\},\\
m_{s_{n_s}\alpha_{n_s}}\, \alpha_{n_s} + \frac{1}{2} m_{s_{n_s} \sigma_{n_s}}. & i=n_s,
\end{array}
\right.
\end{equation}

Therefore, we find that the dual axial RRF equations are simply
\begin{equation}\label{eq:daerrfeqns}
 \boxed{\frac{\partial \alpha_i}{\partial a_i} \, \dot a_i +  \frac{\partial \alpha_i}{\partial a_{i+1}} \, \dot a_{i+1} +  \frac{\partial \alpha_i}{\partial s_i} \, \dot s_i +   \frac{\partial \alpha_i}{\partial s_{i+1}} \, \dot s_{i+1} = -2\, \frac{\epsilon_{s_i}}{s^*_i}.}
 \end{equation}

\subsection{The $\sigma_i$, or dual dodecahedral edge equations}

In many ways, the dual RRF equations for the $\sigma_i$ edges are simpler than the $\alpha_i$ equations because they have no special boundary terms. On the other hand, the Ricci tensors are more complicated. There is one dual dodecahedral edge equation associated each of the $(n_s-1)$ dodecahedron edges $\sigma_i$.  Each of these equations can be expressed locally in terms of the edges $a_i$, $s_i$, and $s_{i+1}$ of the icosahedral frustum model and their time derivatives,
\begin{equation}\label{eq:sigma}
\frac{\dot \sigma_i}{\sigma_i} = -Rc_{\sigma_i}.
\end{equation}
The axial edge, $\sigma_i$, reaches from the circumcenter of one triangle frustum to the adjacent triangle frustum sharing a common trapezoid face sandwiched between the $i$'th  and $(i+1)$'st icosahedrons.  This edge is given by the sum of Eq.~\ref{eq:sigma_1},
\begin{equation}
\sigma_i = \left(\frac{a_i^2 (s_i+s_{i+1})}{\sqrt{3 a_i^2-(s_i-s_{i+1})^2}\sqrt{4 a_i^2-(s_i-s_{i+1})^2}}\right),
\end{equation}
and its time derivative is
\begin{equation}
\dot \sigma_i = \frac{\partial \sigma_i}{\partial a_i} \, \dot a_i +  \frac{\partial \sigma_i}{\partial s_i} \, \dot s_i +   \frac{\partial \sigma_i}{\partial s_{i+1}} \, \dot s_{i+1},
\end{equation}
yielding the left-hand side of Eq.~\ref{eq:sigma}.

The right-hand side of Eq.~\ref{eq:alpha} are the Ricci tensors associated with the dual dodecahedron edges,
\begin{eqnarray}\label{eq:Rcsig1}
Rc_{\sigma_i}  & = & \sum_{\ell_j \in \sigma^*_i} Rc^{(hyb)}_{\ell_j} \left( \frac{V_{\sigma_i \ell_j}}{V_{\sigma_i}} \right) \\ \label{eq:Rcsig2}
& =  & 2 Rc^{(hyb)}_{a_i} \left(  \frac{V_{\sigma_i a_i}}{V_{\sigma_i}} \right) +
Rc^{(hyb)}_{s_i} \left(  \frac{V_{\sigma_i s_i}}{V_{\sigma_i}} \right) +
Rc^{(hyb)}_{s_i} \left(  \frac{V_{\sigma_i s_i}}{V_{\sigma_i}}\right) \\ \label{eq:Rcsig3}
& = &  4\, \frac{\epsilon_{a_i}}{a^*_i} \left(   \frac{V_{\sigma_i a_i}}{V_{\sigma_i}} \right)+
2 \, \frac{\epsilon_{s_i}}{s^*_i}  \left(  \frac{V_{\sigma_i s_i}}{V_{\sigma_i}} \right) +
2 \, \frac{\epsilon_{s_{i+1}}}{s^*_{i+1}} \left(  \frac{V_{\sigma_i s_i}}{V_{\sigma_i}}\right),
\end{eqnarray}
where the ratios of the restricted dual hybrid volumes to the hybrid volumes are expressed in terms of the moment arms in Appendix~\ref{app:1},
\begin{eqnarray}
 \frac{V_{\sigma_i a_i}}{V_{\sigma_i}} & = & \frac{m_{a_i\sigma_i} a_i}{2 m_{a_i\sigma_i} a_i + m_{s_i\sigma_i} s_i + m_{s_{i+1}\sigma_i} s_{i+1}},\\
 \frac{V_{\sigma_i s_i}}{V_{\sigma_i}} & = & \frac{m_{s_i\sigma_i} s_i}{2 m_{a_i\sigma_i} a_i + m_{s_i\sigma_i} s_i + m_{s_{i+1}\sigma_i} s_{i+1}}, \\
 \frac{V_{\sigma_i s_i}}{V_{\sigma_i}} & = &\frac{m_{s_{i+1}\sigma_i} s_{i+1}}{2 m_{a_i\sigma_i} a_i + m_{s_i\sigma_i} s_i + m_{s_{i+1}\sigma_i} s_{i+1}}.
\end{eqnarray}
Two of the three deficit angles in Eq.~\ref{eq:Rcsig3} are given by Eq.~\ref{eq:defs}, while the remaining axial-edge deficit angle  $\epsilon_{a_i}$ is associated with the edge common to five triangle frustum blocks,
\begin{equation}
\epsilon_{a_i} = 2\pi - 5 \theta_{a_i}.
\end{equation}
The dihedral angles ($\theta$) are given in Appendix~\ref{app:1}. Once again, care must be taken when calculating the dual areas, $a^*_i$, as shown in the upper-left of Fig.~\ref{Fig:dualareas}.  This dual area is given by moment arms and the dual edges, $\sigma_i$,  in the boundary of $a^*_i$. It is simply the sum of five identical isosceles triangle areas,
\begin{equation}
a^*_i = \frac{5}{2}\, m_{a_i\sigma_i}\, \sigma_i.
\end{equation}

The dual-dodecahedron-edge RRF equations are
\begin{equation}\label{eq:dderrfeqns}
 \boxed{ \frac{\partial \sigma_i}{\partial a_i} \, \dot a_i +  \frac{\partial \sigma_i}{\partial s_i} \, \dot s_i +   \frac{\partial \sigma_i}{\partial s_{i+1}} \, \dot s_{i+1} = -4\, \frac{\epsilon_{a_i}}{a^*_i} \left(   \frac{V_{\sigma_i a_i}}{V_{\sigma_i}} \right)-
2 \, \frac{\epsilon_{s_i}}{s^*_i}  \left(  \frac{V_{\sigma_i s_i}}{V_{\sigma_i}} \right) -
2 \, \frac{\epsilon_{s_{i+1}}}{s^*_{i+1}} \left(  \frac{V_{\sigma_i s_i}}{V_{\sigma_i}}\right).}
 \end{equation}

\section{the Numerical Algorithm and Simulations}\label{sec:numerics}

We report on our numerical solution of Eqs.~\ref{eq:daerrfeqns} and \ref{eq:dderrfeqns} given the initial profile dual-loabed profile of Eq.~\ref{Eq:ivd}.  These equations form a sparsely-coupled first-order system of nonlinear algebraic equations,
\begin{equation}\label{eq:matrixeq}
M_{ij} \, \dot \ell_j = f_i,\ \ \ \forall i,j\in\{1,2,\ldots,2n_s-1\}.
\end{equation}
We define the coordinate vector as alternating icosahedral and axial edges,
\begin{equation}
\pmb{\ell} = \{s_1,a_1,s_2,a_2, \ldots, s_{N_{s-1}},a_{N_{s-1}},s_{n_s}\};
\end{equation}
consequently, the square matrix $\pmb M$ is sparse with four non-zero diagonals and populated with partial derivatives of the dual edges with respect to isosceles frustum lattice edges. The vector $\pmb f$'s components are the weighted sectional curvatures.  We solve the equations using a fixed-time step 4-th order Runge-Kutta algorithm and remesh the lattice occasionally to keep the circumcenters roughly inside the frustum blocks.  As a diagnostic, we linearize,
\begin{equation} \label{eq:jac}
{\pmb J} := \delta\left({\pmb M}^{-1} \cdot {\pmb f}\right),
\end{equation}
and track this Jacobian and its corresponding eigenvalues during the evolution.

The initial radial profile of the double-lobed geometry given by Eq.~\ref{Eq:ivd} is translated into initial values for the axial and icosahedral edges using Eqs.~\ref{eq:s2r} and \ref{eq:a2ac}. However, we exploit our freedom in the axial placement of each icosahedron.  We choose to concentrate more icosahedra near the waist where we expect a singularity in finite time. In order to remesh, we construct an cubic interpolating function, $s\left( a^{coord}\right)$  for the simplicial geometry.  We also extrapolate to the poles ($s=0$) of the dumbbell geometry in order to keep the end cap icosahedra bounded away from the poles. This extended interpolating function allows us to place the icosahedra along the geometry based on a distribution function.  We choose a Gaussian distribution centered on the waist,
\begin{equation}
a^{(new)}_i =  a^{(ext)}_i
\left(
1 - \kappa \exp{
\left( -\frac{  (\bar a^{{coord}}_i /a^{(coord)}_{max})^{2}  }{  2 \sigma^{2}  } \right)
}
\right)
\end{equation}
Here we use $\kappa=0.95$, $\sigma=0.1$, $a^{(coord)}_{max}$ is the maximal length geodesic from pole to pole,  $\bar a^{{coord}}_i  = \left( a^{{coord}}_{i+1} - a^{{coord}}_i   \right)/2$, and $a^{(ext)}_i = a^{(coord)}_{max}/(n_s+1)$ represents the equally-spaced $a$'s along the extrapolated length of the interpolating function.

During the evolution, we want each of the circumcenters to lie within, or nearly within, their respective triangular frustum blocks.  The condition for this well-centeredness is obtained from Eqs.~\ref{eq:alpha_i}, \ref{eq:alpha_ip}, and an expression for the altitude of the frustum block; we impose,
\begin{equation}\label{eq:wc}
a_i \ge \sqrt{\frac{s_{i+1} |s_i-s_{i+1}|}{3}}.
\end{equation}
Given the very poor azimuthal resolution afforded to us by the icosahedra, this well-centeredness condition is nearly impossible to satisfy, and we have to allow each circumcenter to evolve slightly outside its frustum block before we remesh.  We find, on occasion, the adjacent circumcenters cross in very short time scales, and this abruptly crashes the evolution. This problem is exacerbated when we increase the number of icosahedral cross sections in our model. In this case the height of the frustum block vanishes, we observe that the gradients in the radial profile, $(s_i+1-s_i)/a_i$,  will cause changes in the axial edges, $a_i$, which in turn cause the circumcenter to evolve through the largest equilateral triangle face.  With this in mind, we find empirically that the nearly maximal number of icosahedral cross sections we can use in our simulation in order to to evolve stably is  $n_s=45$. However, if we were to remove the condition that the corrections be limited to regular icosahedron then the azimuthal resolution could increase which would allow cross-sectional slices, which should lead more finely resolved evolution and comparison with the continuum. However, this was beyond the scope and intent of this paper.  

\begin{figure}
 \includegraphics[width=14cm]{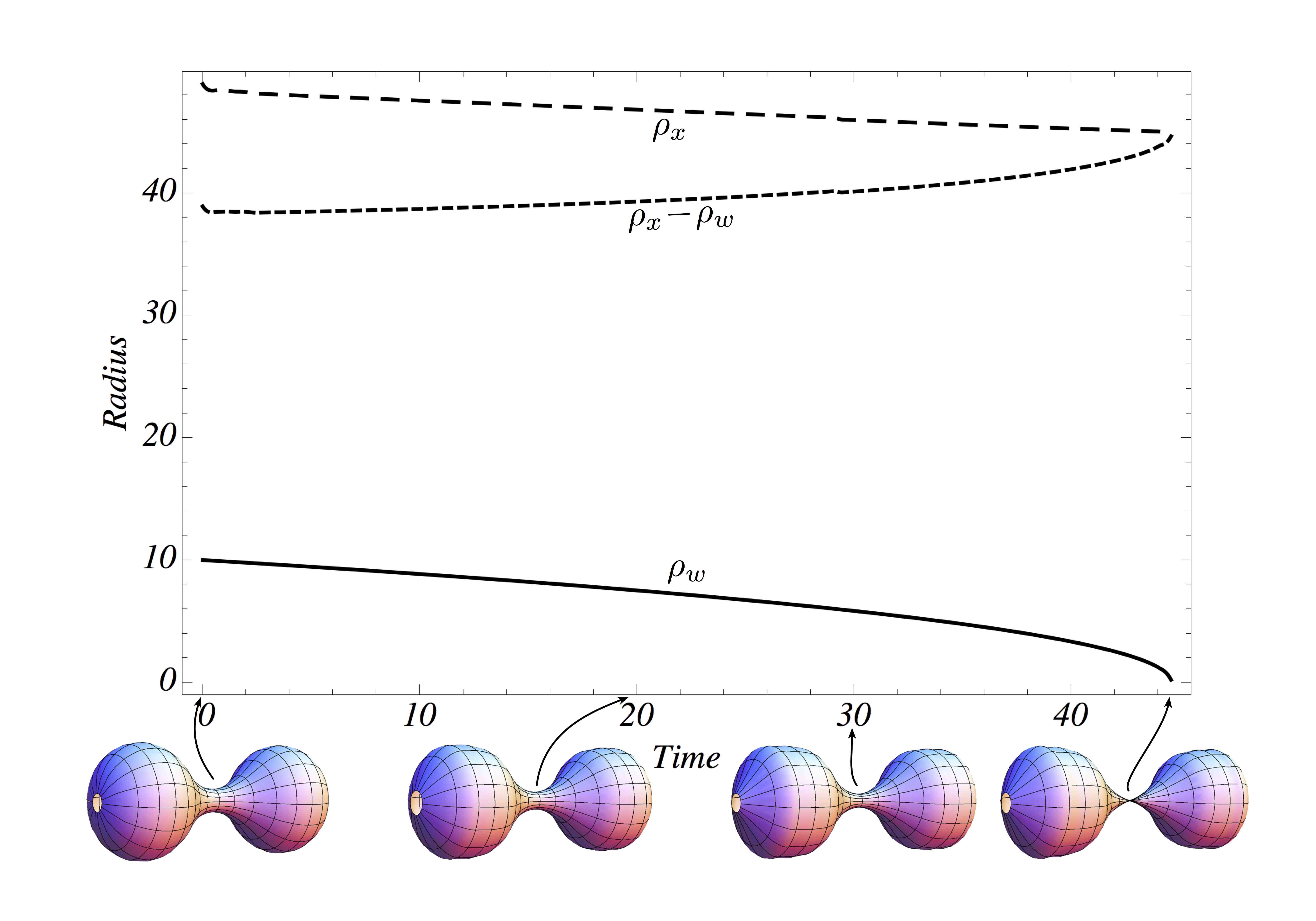}
\caption{The evolution of the 3-dimensional SRF icosahedral frustum geometry using the RRF equations with the remshing described above. The solid line is the radius of the waist as a function of time. The long dashed line at the top is the time evolution of the largest radius of the lobes. The short dashed line in between the other two lines is the difference between the maximum radius and the waist radius, thus demonstrating a neck pinching singularity.  Just below the time axis, and at four separate times during the evolution we plot an embedding surface representing the dumbbell geometry as it evolves toward neck pinch singularity.  }
\label{fig:p75}
\end{figure}

In Fig.~\ref{fig:p75} we show the results of our evolution.  This evolution includes remeshing at those times during the evolution when the well-centeredness condition (Eq.~\ref{eq:wc}) is violated.  It is interesting to note that we observed that this evolution only runs a few units in time without remeshing.  Furthermore, we find that the system of dual-edge RRF equations is rather stiff, with the condition number of the Jacobian in Eq.~\ref{eq:jac} becomes on the order of ${\cal C} \sim 10^8$.  As a rule of thumb, the accuracy is diminished by $\log_{10}\left({\cal C}\right)$.  Therefore, it is important to keep the condition number as small as possible.  Remeshing reduces the condition number by a factor of 2 and also produces more negative eigenvalues in the Jacobian.  While the equations we are solving are obtained by the method of lines, the equations are nonlinear. The eignenvalue spectrum of the Jacobian will change during evolution; that we observe more negative eigenvalues under remeshing  is a positive diagnostic.  We identify an interesting phenomenon where the axial distance of every other icosahedron shrinks, while the adjacent axial edges grow. The relatively poor resolution of the icosahedra prohibits us from exploring this in much detail, and a refinement in the azimuthal directions is beyond the scope and intent of this paper. 

\begin{figure}
 \includegraphics[width=14.5cm]{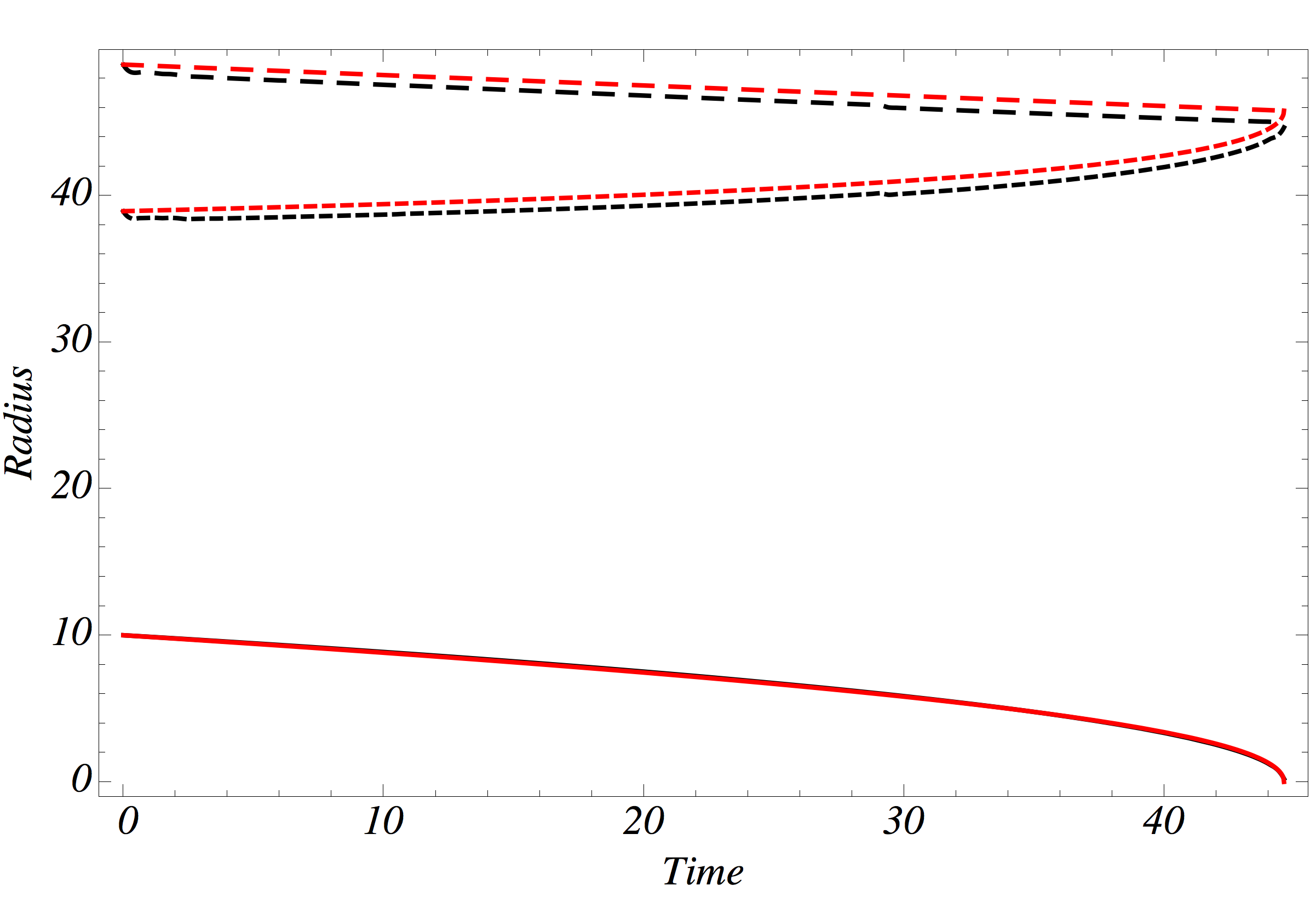}
\caption{We show here a comparison of the continuum RF equation and the RRF equations. The black curves are the solution of the RRF equations, and the red curves are as solution the continuum equations.  The solid lines (lying on top of each other) are the radius of the waist as a function of time. The long dashed lines at the top are the time evolution of the largest radius of the lobes.  The short dashed lines in between are the difference between the maximum radius and the waist radius, thus demonstrating a neck pinch singularity.  These curves are only scaled in time so that the pinch happened at the same instant. Given the relatively low resolution afforded by the icosahedral cross sections, the agreement is rather striking.  }
\label{fig:p75Steps}
\end{figure}

We also solve the continuum RF equation, Eq.~\ref{eq:rhoeqn}, for the same initial radial profile and compare it to the SRF solution. These two solutions show excellent agreement and demonstrate neck pinching in SRF.  To make the comparison shown in Fig.~\ref{fig:p75Steps}, we only rescale the collapse time of the continuum equations to agree with the SRF evolution.

\section{ Toward Higher Resolution, Proper  Sampling and Surgery in  SRF}\label{sec:conclusion}

The purpose of this paper has been to demonstrate a neck pinching evolution in 3-dimensional SRF.  Our comparison with the continuum solution accomplishes this goal. We have learned that these equations require delicate numerical techniques and higher resolution. We have demonstrated here the necessity of remeshing, and we expect to need a more sophisticated adaptive mesh and remeshing algorithm for higher-dimensional SRF simulations.  Our model used regular icosahedral cross sections and therefore we could not explore the scaling behavior of our code with resolution in a meaningful way. However, the extension of our warped-produce geometries to higher dimensions will not change the structure or complexity of the RRF equations.  There will still be only one axial equation and one azimuthial equation per cross-sectional slice.   It  should not add any substantial computational demands to extend beyond 3-dimensions.  The sparsity of the matrix we need to invert will not increase in the discretization of our warped-product geometry. Notwithstanding this observation, the availability of regular polytopes in four and higher dimensions is restricted to the cross polytope. The azimuthal resolution would be more course than our 3-dimensional icosahedral model and we would not be able to resolve the dumbbell along the axial direction as well either. The computational scaling of our algorithm should be dominated by the inversion of the a sparse square matrix in Eq.~\ref{eq:matrixeq}.  This matrix needs to be constructed and inverted at each time step. The structure of this matrix, $M_{i j}$,  is determined by the lattice structure. In our case it will have at most four non-zero diagonals.  This is also the case for higher dimensional dumbbell geometries.  A regular lattice structure will yield a well structured matrix. There are efficient parallelized preconditioned matrix inversion algorithms to extend our numerical techniques for arbitrary geometries in three and higher dimensions \cite{Chow:2001}.  We see no way in three and higher dimensions to decouple the RRF equations.

Two conditions govern the proper sampling of the geometry.  First, we would like to have large blocks where the curvature is small and small blocks where the curvature is large. This can be accomplished with the condition that all the deficit angles be small and roughly equal, since the sectional curvature will then be inversely proportional to its dual area. Second, we would like to keep the blocks as well-centered and equilateral as possible.  This can be accomplished by using a well-centeredness or fullness condition introduced by Whitney \cite{Whitney:1957}, and others used in RC, e.g. the waste function \cite{Miller:1986}.

We consider that the best way to advance this numerical work in SRF follows three distinct approaches. First, we can introduce a higher-resolution triangulation of each of the cross-sectional spheres in the current model. We are currently exploring a simplicial model that can be arbitrarily refined, and refer to this as the ``continuum SRF model" \cite{Miller:2014}. With this higher azimuthal resolution, we will be better able to maintain well-centeredness of the lattice to a higher tolerance, and we hope to introduce more spherical cross-sectional spheres.  The numerical intricacies in solving the RRF equations that we may learn from the continuum model will guide us to a fully-generic (e.g., no axial symmetry restrictions) simplicial geometry evolving under SRF. Note that for a general simplicial geometry there are more dual edges than simplicial edges which necessitates solving the simplicial-edge RRF equations directly. In addition, we need to develop a sophisticated implementation of remeshing as well as adapting the mesh during evolution (i.e. resample the mesh in areas at higher or lower resolution).  We also need an adaptive time step stiff integrator.  This example makes clear that such techniques are essential for the analysis of singularity development in this model. Such techniques will enable us to automatically detect singularity formation, perform surgery in these regions, and continue our integration. Some of us are already exploring surgery techniques in SRF \cite{Corne:2014}.   Finally, we think it important to explore a reformulation of SRF where we have a diffeomorphically equivalent flow under a convex energy functional, e.g. a simplicial application of the de~Turck trick, or a simplicial version of a convex entropy functional \cite{LiYau:1986,Perelman:2003}.

\section*{Ackowledgements}

We thank Dan Knopf, David Glickenstein, Huai-Dong Cao, Howard Blair, Chris Tison,  and Arkady Kheyfets for discussions, and Xuping Wang for his suggestions on the SRF equations for this model.  We are especially grateful for the discussions we had with Jonathan McDonald.  This work has benefited from the discussions and input we received from Lars Hernquist on mesh refinement and re-meshing strategies, some of which we partially implemented here. WAM would like to thank the Department of Mathematics at Harvard University for their support and hospitality this year.  This work was supported from USAF Grant \# FA8750-11-2-0089.  WAM, CT, and SR acknowledge support from Air Force Office of Scientific Research through the American Society for Engineering Education's 2012 Summer Faculty Fellowship Program, and from AFRL/RITA and the Griffiss Institute's 2013 Visiting Faculty Research Program.  Any opinions, findings and conclusions or recommendations expressed in this material are those of the author(s) and do not necessarily reflect the views of the AFRL.

\begin{appendix}
\section{The Geometry of the Equilaterial Triangle Frustum Block}\label{app:1}

\begin{figure}
 \includegraphics[width=7.5cm]{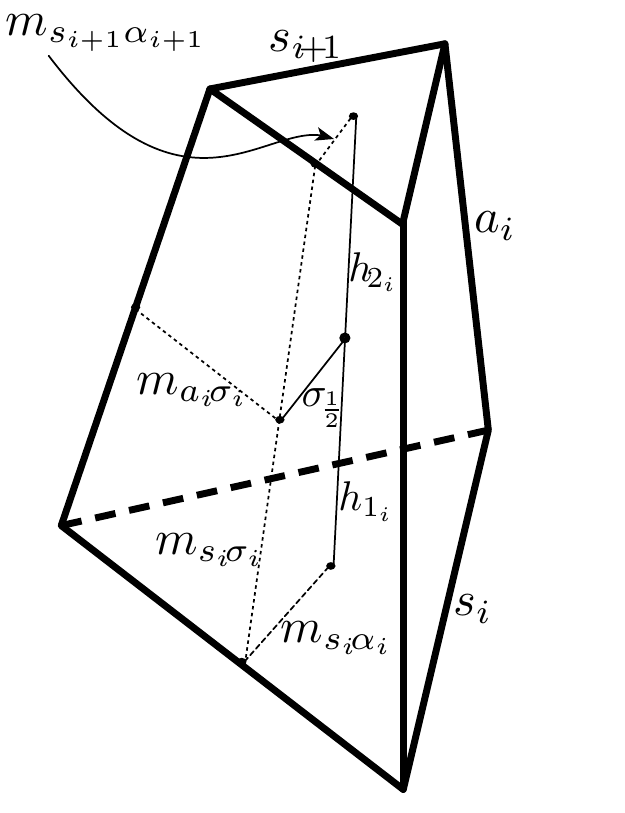}
\caption{The isosceles triangle frustum block}
\label{Fig:frustum}
\end{figure}

The SRF equations depend directly on the geometry of the triangular frustum block. We will outline the relevant geometric features of this polyhedron that are used to construct the axial and icosahedral-edge RRF equations for this model.  We focus in this section on a single isosceles-based triangular frustum block as illustrated in Fig.~\ref{Fig:frustum}.
Here the block has three types of edges:  three axial edges each of length $a_i$, three equal edges of the base isosceles triangle of length $s_i$, and three equal edges of the cap triangle of length $s_{i+1}$.  Here we assume the base triangle is larger than the top cap triangle, $s_i > s_{i+1}$.  Since all three of the axial edges are equal then the top triangle is parallel to the base triangle; also, the circumcenter of the two triangles and the circumcenter of the frustum block all lie on the same line.

Each of our SRF equations are constructed, in part, from the three dihedral angles of the frustum block.  The dihedral angle between any two of its three trapezoidal faces sharing an edge, $a$, is the axial dihedral angle,
\begin{equation}\label{eq:theta_a}
\theta_{a_i} = \arccos{\left(\frac{2 a_i^2-(s_i-s_{i+1})^2}{4 a_i^2-(s_i-s_{i+1})^2}\right)}.
\end{equation}
The dihedral angle along edge $s_i$ is the angle between the base triangle and the trapezoid sharing $s_i$ we refer to as the base dihedral angle,
\begin{equation}\label{eq:theta_{s_i}}
\theta_{s_i} = \arccos{\left(\frac{\sqrt{3} (s_i-s_{i+1})}{3\sqrt{4 a_i^2-(s_i-s_{i+1})^2}}\right)},
\end{equation}
and consequently the corresponding dihedral angle associated to the top cap edge, $s_{i+1}$, is the supplementary angle,
\begin{equation}\label{eq:theta_{s_{i+1}}}
\theta_{s_{i+1}}=\pi-\theta_{s_i}.
\end{equation}

This frustum block contains the segments of three of the dual edges in our model:  the dual decahedron edge, $\sigma_i$, the dual axial edge $\alpha_i$, and the dual axial edge $\alpha_{i+1}$.  The circumcentric dual edge $\sigma_i$ associated with this frustum block is the line segment that starts at  the circumcenter of this frustum block, pierces through the circumcenter of one of the three trapezoidal blocks of the frustum, and terminates at the circumcenter of the adjacent frustum block.  In this case half of the dual edge, $\sigma_{\frac{1}{2}_i}$, lies in the frustum considered here,
\begin{equation}\label{eq:sigma_1}
\sigma_{\frac{1}{2}_i} = \frac{1}{2} \left(\frac{a_i^2 (s_i+s_{i+1})}{\sqrt{3 a_i^2-(s_i-s_{i+1})^2}\sqrt{4 a_i^2-(s_i-s_{i+1})^2}}\right),
\end{equation}
and $\sigma_i=2\, \sigma_{\frac{1}{2}_i}$. There are two segments of dual axial edges within the frustum, one dual to the base isosceles triangle, $\alpha_{i}$, and the other with the top cap triangle, $\alpha_{i+1}$. Here, the dual axial edge segment reaches from the circumcenter of the frustum to the circumcenter of the base isosceles triangle,
\begin{equation}\label{eq:alpha_i}
h_{1_i} = \frac{\sqrt{3}}{6} \left( \frac{3 a_i^2-2 s_i (s_i-s_{i+1})}{\sqrt{3 a_i^2-(s_i-s_{i+1})^2}}\right).
\end{equation}
The dual edge reaching from the circumcenter of the frustum to the circumcenter of the top cap isosceles triangle can be obtained from $\alpha_{i+1}$ by switching the base edges in Eq.~\ref{eq:alpha_i}, i.e. $s_{i+1} \leftrightarrow s_i$,
\begin{equation}\label{eq:alpha_ip}
h_{2_{i+1}} = \frac{\sqrt{3}}{6} \left( \frac{3 a_i^2+2 s_{i+1} (s_i-s_{i+1})}{\sqrt{3 a_i^2-(s_i-s_{i+1})^2}}\right).
\end{equation}

Each of our SRF equations depend on one or more of the five types of moment arms. A moment arm, e.g. $m_{a_i\sigma_i}$ is defined as the line segment reaching from the center of an edge of the frustum, $a_i$ in this case, to the perpendicular bisector of a dual edge, $\sigma_i$ in this case. Since $\sigma_i$ is dual to a trapezoid, there are three moment arms, $m_{a_i \sigma_i}$, $m_{s_i \sigma_i}$ and $m_{s_{i+1} \sigma}$. Moment arm $m_{a_i \sigma_i}$ is the length of the edge that  reaches from the center of edge, $a_i$ to the circumcenter of the trapezoidal face of the frustum,
\begin{equation}\label{eq:m_asigma}
m_{a_i\sigma_i} = \frac{a_i (s_i+s_{i+1})}{2 \sqrt{4 a_i^2-(s_i-s_{i+1})^2}}.
\end{equation}
This moment arm is perpendicular to dual edge $\sigma_i$.  Similarly,  $m_{s_i \sigma_i}$ is the length of the edge that  reaches from the center of either edge $s_i$ to the circumcenter of the trapezoidal face of the frustum,
\begin{equation}\label{eq:m_bsigma}
m_{s_i\sigma_i} = \frac{2 a_i^2 -s_i (s_i-s_{i+1})}{2 \sqrt{4 a_i^2-(s_i-s_{i+1})^2}},
\end{equation}
and finally,
\begin{equation}\label{eq:m_csigma}
m_{s_{i+1}\sigma_i} = \frac{2 a_i^2 +s_{i+1} (s_i-s_{i+1})}{2 \sqrt{4 a_i^2-(s_i-s_{i+1})^2}},
\end{equation}
which is obtained from the base moment arm via  $s_{i+1} \leftrightarrow s_i$.  The two remaining moment arms are associated with the dual axial edges, $\alpha_i$.  In particular, the moment arm $m_{s_i\alpha_i}$ reaches from the center of edge $s_i$ to the circumcenter of the base triangle.  It is the inradius of the base triangle,
\begin{equation}\label{eq:m_balpha}
m_{s_i\alpha_i} = \frac{\sqrt{3}}{6}\, s_i.
\end{equation}
The corresponding moment arm associated with the top cap of the frustum is also the inradius of the top triangle,
\begin{equation}\label{eq:m_calpha}
m_{s_{i+1}\alpha_{i+1}} = \frac{\sqrt{3}}{6}\, s_{i+1}.
\end{equation}

The only other quantities that enter into the normalized SRF equations from the frustum block is its volume, $V_f$,
\begin{equation}\label{eq:volume_f}
V_{f_i} =  \frac{1}{12} \left(s_i^2+s_{i+1}^2 + s_i s_{i+1}\right)\, \sqrt{3 a_i^2-(s_i-s_{i+1})^2}.
\end{equation}
In addition, our frustum geometry has two regular icosahedral boundaries at both ends of the dumbbell.   The dihedral angle of the regular icosahedron of edge length, $s$,
\begin{equation}\label{eq:theta_icosa}
\theta_{icosa} =  \arccos{ \left(-\frac{\sqrt{5}}{3}\right)},
\end{equation}
as well as its 3-volume,
\begin{equation}\label{V_icosa}
V_{icosa} = \frac{5}{12} \left( 3+\sqrt{5} \right) s^3,
\end{equation}
also enter into our axial and spherical SRF equations.

\section{The Equivalence of the RRF and Dual-Edge RRF Equations for Axisymmetric Geometries} \label{app:B}

In this section we prove that the RRF equations associated with the edges of the icosahedral-frustum model  are equivalent to the dual-dodecahedron-edge RRF equations.  We follow as closely as possible the notation and equations  from recent work in Sec.~3 and Sec.~5 of the Simplicial Ricci Flow manuscript \cite{Miller:2013}.
\vskip 0.15 true in
\noindent
{\em Theorem:}  The Regge-Ricci Flow equations associated to the edges of the icosahedral-frustum-block geometry are equivalent to the dual-edge Regge-Ricci Flow equations.
\vskip 0.1 true in
\noindent
{\em Proof.}\newline \noindent
 There are two sets of RRF equations associated with the edges of the icosahedral-frustum model as shown in the right-hand side of Eq.~\ref{eq:drrf-rrf}. One equation, the $s_i$-equation,  is associated with each icosahedron, the other set is associated with the axial edges ($a_i$). In particular, there is one $a_i$-equation associated with each axial edge. These can be expressed in terms of the moment arms, deficit angles, dual areas and edge lengths. All moment arms are strictly positive.  The $s_i$-equation can conveniently be written as  the sum of three terms,
\begin{align}\label{eq:si}
  m_{s_i\sigma_i}  &  \underbrace{\left( \dot \sigma_{i} + \left( \frac{\sigma_{i}}{\sigma^*_{i}} \right) \left\{ \left( \frac{s_{i}}{s^*_{i}} \right) m_{s_{i}\sigma_{i}} \epsilon_{s_{i}} + \left( \frac{s_{i+1}}{s^*_{i+1}} \right) m_{s_{i+1}\sigma_{i}}\epsilon_{s_{i+1}} + 2 \left( \frac{a_{i}}{a^*_{i}} \right) m_{a_{i}\sigma_{i}}\epsilon_{a_{i}} \right\} \right)}_{Term\ \sigma_i} \nonumber \\
+  m_{s_i\sigma_{i+1}} & \underbrace{\left( \dot \sigma_{i+1} + \left( \frac{\sigma_{i+1}}{\sigma^*_{i+1}} \right) \left\{ \left( \frac{s_{i+1}}{s^*_{i+1}} \right) m_{s_{i+1}\sigma_{i+1}} \epsilon_{s_{i+1}} + \left( \frac{s_{i+2}}{s^*_{i+2}} \right) m_{s_{i+2}\sigma_{i+1}}\epsilon_{s_{i+2}} + 2 \left( \frac{a_{i+1}}{a^*_{i+1}} \right) m_{a_{i+1}\sigma_{i+1}}\epsilon_{a_{i+1}} \right\} \right)}_{Term\ \sigma_{i+1}} \nonumber \\
+  2 m_{s_i\alpha_i} & \underbrace{\left( \dot \alpha_i + 3 \, \left(\frac{\alpha_i}{\alpha^*_i} \right)  \left( \frac{s_i}{s^*_i} \right)m_{s_i\alpha_i} \epsilon_{s_i} \right)}_{Term\ \alpha_i} = 0.
\end{align}
The axial-edge equation associated with edge $a_i$ is just {\em Term}~$\sigma_i$ of the $s_i$-equation.
\begin{equation}\label{eq:ai}
\underbrace{\left( \dot \sigma_{i} + \left( \frac{\sigma_{i}}{\sigma^*_{i}} \right) \left\{ \left( \frac{s_{i}}{s^*_{i}} \right) m_{s_{i}\sigma_{i}} \epsilon_{s_{i}} + \left( \frac{s_{i+1}}{s^*_{i+1}} \right) m_{s_{i+1}\sigma_{i}}\epsilon_{s_{i+1}} + 2 \left( \frac{a_{i}}{a^*_{i}} \right) m_{a_{i}\sigma_{i}}\epsilon_{a_{i}} \right\} \right)}_{Term\ \sigma_i} = 0
\end{equation}
The under braces on these equations are important for the clarity of the proof.
\newline \noindent
There are just two dual-dodecahedral-edge RRF equations for this model.  The $\alpha_i$-equation is simply {\em Term} ~$\alpha_i$,
\begin{equation}\label{eq:alpi}
\underbrace{\left( \dot \alpha_i + 3 \, \left(\frac{\alpha_i}{\alpha^*_i} \right)  \left( \frac{s_i}{s^*_i} \right)m_{s_i\alpha_i} \epsilon_{s_i} \right)}_{Term\ \alpha_i} = 0,
\end{equation}
and the $\sigma_i$-equation is simply {\em Term} ~$\sigma_i$,
\begin{equation}\label{eq:sigi}
\underbrace{\left( \dot \sigma_{i} + \left( \frac{\sigma_{i}}{\sigma^*_{i}} \right) \left\{ \left( \frac{s_{i}}{s^*_{i}} \right) m_{s_{i}\sigma_{i}} \epsilon_{s_{i}} + \left( \frac{s_{i+1}}{s^*_{i+1}} \right) m_{s_{i+1}\sigma_{i}}\epsilon_{s_{i+1}} + 2 \left( \frac{a_{i}}{a^*_{i}} \right) m_{a_{i}\sigma_{i}}\epsilon_{a_{i}} \right\} \right)}_{Term\ \sigma_i} = 0.
\end{equation}
\newline\noindent
Suppose Eqs.~\ref{eq:si} and \ref{eq:ai} are satisfied.  Then, Eq. \ref{eq:sigi} is automatically satisfied as it equals Eq. \ref{eq:ai}.  Also, by an indexing argument, Eq. \ref{eq:ai} implies that \emph{Term} $\sigma_{i} = 0 $ and \emph{Term} $\sigma_{i+1} = 0 $.  The vanishing of \emph{Term} $\sigma_{i} $ and \emph{Term} $\sigma_{i+1} $, with Eq. \ref{eq:si} and the strict positivity of moment arms, imply that \ref{eq:alpi} is true.
\newline \noindent
Conversely, if the dual-edge equations (Eqs.~\ref{eq:alpi} and \ref{eq:sigi}) are satisfied, then all three terms in Eq. \ref{eq:si} vanish so that their sum is zero.  Also, Eq.~\ref{eq:ai} is identical to Eq. \ref{eq:sigi}.
\newline\noindent
QED
\end{appendix}


\begin{thebibliography}{99}
\bibitem{Hamilton:1982} R. Hamilton, ``Three-manifolds with positive Ricci curvature,'' {\em J. Diff. Geom} {\bf 17} (1982), 255-306.
\bibitem{Cao:2003}H-D. Cao, B. Chow, S-C Chu \& S-T Yau, eds., {\em Collected Papers on Ricci Flow}  in Series in Geometry and Topology, Volume 37 (International Press; Somerville, MA; 2003).
\bibitem{Chow:2004} B. Chow \& D. Knopf, {\em The Ricci Flow: An Introduction}, Mathematical Surveys and Monographs, Volume 110 (American Mathematical Society; Providence, RI; 2004).
\bibitem{Chow:2006} B Chow, P. Lu \& L. Ni, {\em Hamilton's Ricci Flow}, Graduate Studies in Mathematics, Volume 77 (American Mathematical Society; Providence, RI; 2006).
\bibitem{Chow:2007} B. Chow, S-C Chu, D. Glickenstein, C. Guenther, J. Isenberg, T. Ivey, D. Knopf, P. Lu, F. Luo \& L. Ni, {\em The Ricci Flow: Techniques and Applications, Part 1: Geometric Aspects}, Mathematical Surveys and Monographs, Volume 135 (American Mathematical Society; Providence, RI; 2007).
\bibitem{Gu:2012} X. Yu, X. Yin, W. Han, J. Gao \& X. Gu, ``Scalable routing in 3D high genus sensor networks using graph embedding,'' {\em  INFOCOM 2012}: 2681-2685;
Y. Wang, J.  Shi, X.  Yin, X.  Gu, T. F. Chan, S-T Yau, A. W. Toga \& P. M. Thompson, ``Brain surface conformal parameterization with the Ricci flow,'' {\em IEEE Trans. Med. Imaging} {\bf 31}(2)  (2012) 251-264.
X. Gu, F.  Luo \& S-T Yau, ``Fundamentals of computational conformal geometry,'' {\em Mathematics in Computer Science} {\bf 4}(4) (2010) 389-429;
B. Chow \& F. Luo, ``Combinatorial Ricci flows on surfaces,'' {\em J. Differential Geometry} {\bf 63} (2003) 97-129.
\bibitem{Peiro:2005} J. Peiro \& S. Sherwin, {\em Finite Difference, Finite Element and Finite Volume Methods For Partial Differential Equations,}  in Handbook of Materials Modeling, Volume 1,  Methods and Models, Springer, 2005.
\bibitem{Humphries:1997} S. Humphries, Jr., {\em Finite-Element Methods for Electromagnetism}, http://www.fieldp.com/freeware/finite\_element \_electromagnetic.pdf; originally published as {\em Field Solutions on Computers}  (ISBN 0-8493-1668-5) (Taylor and Francis, Boca Raton, 1997).
\bibitem{Regge:1961} T. Regge, ``General relativity without coordinates,'' {\em Il Nuovo Cimento} {\bf 19} (1961) 558-571.
\bibitem{Gentle:1998} A. P. Gentle \& W. A. Miller, ``A fully (3+1)-D Regge calculus model of the Kasner cosmology,'' {\em Class.Quant.Grav.} {\bf 15} (1998) 389-405.
\bibitem{Hirani:2005} M. Desbrun, A. N. Hirani, M. Leok \&  J. E. Marsden, ``Discrete exterior calculus,''
e-print arXiv:math/0508341v2 [math.DG] on arxiv.org (2005).
\bibitem{Glickenstein:2011a} D. Glickenstein, D. Champion and A. Young, ``Regge's Einstein-Hilbert functional on the double tetrahedron,''  {\em Differential Geom. Appl.} { \bf 29} (2011), 109-124, doi:10.1016/jdifgeo.2010.10.001.
\bibitem{Glickenstein:2011} D. Glickenstein, ``Discrete conformal variations and scalar curvature on piecewise flat two- and three-dimensional manifolds," {\em J. Diff. Geom.} {\bf 87} (2011) 201-238.
\bibitem{G:2005} D. Glickenstein, ``Geometric triangulations and discrete Laplacians on manifolds,'' arXiv:math/0508188 [math.MG].
\bibitem{Ge:2013} H. Ge, ``Discrete Quasi-Einstein Metrics and Combinatorial Curvature Flows in 3-Dimension,'' 	 arXiv:1301.3398 [math.DG].
\bibitem{Forman:2003}R. Forman, ``Bochner's Method for Cell Complexes and Combinatorial Ricci Curvature,'' {\em Discrete Comput. Geom.} {\bf 29} (2003) 323Ð374.
\bibitem{Glickenstein:2005} D. Glickenstein, ``A combinatorial Yamabe flow in three dimensions,'' {\em Topology} {\bf 44} (2005) pp.~791-808.
\bibitem{GuSaucan:2013} D. X. Gu \& E. Saucan, ``Metric Ricci curvature for PL manifolds,'' {\em Geometry} (2013) 694169.
\bibitem{Besson:2007}  G. Besson, ``The geometrization conjecture after R. Hamilton and G. Perelman,'' {\em Rend. Semin. Mat. Univ. Politec. Torino}{\bf 65} (2007) p. 397?411.
\bibitem{Corne:2014} M. Corne, P. Alsing, H. Blair, G. Jones, W. A. Miller, \& V. Nanda, ``Applications of Persistent Homology to Ricci Flow on $S^2$ and $S^3$. 
\bibitem{Mischaikow:2013} Private communication K. Mischaikow, V. Nanda \& J. Bush (2013).
\bibitem{Miller:2013} W. A. Miller, J. R. McDonald, P. M. Alsing, D. Gu \& S-T Yau, ``Simplicial Ricci Flow,''  submitted to {\em Comm. Math. Phys.} (2013); arXiv:1302.0804v1 [math.DG].
\bibitem{AMM:2011} P. M. Alsing, J. R. McDonald \& W. A. Miller, ``The simplicial Ricci tensor," {\em Class. Quantum Grav.} {\bf 28} (2011) 155007 (17 pp).
\bibitem{McDonald:2012} J. R. McDonald, W. A. Miller, P. M Alsing, X. D. Gu, X. Wang \& S-T Yau, ``On exterior calculus and curvature in piecewise-flat manifolds,'' {\em paper submitted to J. Math. Phys.} (2012) arxiv.org/abs/1212.0919.
\bibitem{Knopf:2004} S.  Angenent \& D. Knopf, ``An example of neckpinching for Ricci flow on $S^{n+1}$'' .  {\em Math. Res. Lett.} {\bf 11} (2004) 493-518.
\bibitem{Simon:2000} M. Simon,  ``A class of Riemannian manifolds that pinch when evolved by Ricci flow,'' {\em Manuscripta Math.} {\bf 101}, no. 1  (2000)
89Ð114.
\bibitem{Knopf:2011} S. Angenent, J. Isenberg and D. Knopf, ``Formal matched asymptotics for degenerate Ricci flow neckpinches,'' {\em Nonlinearity} {\bf 24} (2011), 2265-2280.
\bibitem{Cao:2011} H.-D. Cao, {\em Private communication,} UBC, Vancouver BC (2011).
\bibitem{Hamilton:1995} R. S. Hamilton, ``The formation of singularities in the Ricci flow,'' {\em Surveys in Differential Geometry} {\bf 2} (1995) 7-136.
\bibitem{Chow:2003} B. Chow and F. Luo, ``Combinatorial Ricci Flows on Surfaces,'' {\em J. Differential Geom.}  {\bf 63}, no.~1 (2003) 97-129.
\bibitem{LinYau:2010} Y. Lin and S-T Yau, ``Ricci curvature and eigenvalue estimate on locally finite graphs,'' {\em Math. Res. Lett.} {\bf 17} (2010) 343Ð356.
\bibitem{Knopf:2009} D. Knopf, ``Estimating the trace-free Ricci tensor in Ricci flow,''  {\em Journal: Proc. Amer. Math. Soc.} {\bf 137} (2009), 3099-3103.
\bibitem{Thurston:1997} W.  Thurston, {\em Three-dimensional geometry and topology,} Vol. 1. Edited by Silvio Levy, Princeton Mathematical Series, 35,  (Princeton University Press, Princeton, NJ, 1997).
\bibitem{Whitney:1957} H. Whitney, {\em Geometric Integration Theory} (Princeton University Press, Princeton, NJ; 1957) pp. 124Ð135.
\bibitem{Miller:1986} W. A. Miller, ÒGeometric Computation: Null-Strut Geometrodynamics and the Inchworm Algorithm,Ó in {\em Dynamical Spacetimes and Numerical Relativity}, ed. J. Centrella, Cambridge Univ. Press.(1986) 256-303.
\bibitem{Miller:2014} W. A. Miller,  P. M. Alsing, M. Corne \& S. Ray ``An Analysis of Neck Pinching in High Resolution Simplicial Ricci Flows,''  {\em in preparation} (2014).
\bibitem{Corne:2014} M. Corne, P. Alsing and W. Miller, ``Surgery in Simplicial Ricci Flow," {\em in preparation} (2014).
\bibitem{LiYau:1986} P. Li and S-T Yau, ``On the parabolic kernel of the Schr\"odinger operator," {\em Acta Math.} {\bf 156}, no. 3-4 (1986)153Ð201.
\bibitem{Perelman:2003} G. Perelman, ``The entropy formula for the Ricci flow and its geometric applications,'' preprint,
math.DG/0211159; G. Perelman, ``Ricci flow with surgery on three-manifolds,'' preprint, math.DG/0303109; \&
G. Perelman, ``Finite extinction time for the solutions to the Ricci flow on certain three-manifolds,''  preprint, math.DG/0307245.
\bibitem{Chow:2001} R. Chow, ``Parallel implementation and practical use of sparse approximate inverse pre conditioners with a priori sparsity patterns,'' {\em Intl. J. High Perf. Comput. Appl.} {\bf 15} (2001) pp. 56-74.

\end{thebibliography}
\end{document}